# KREIN'S SPECTRAL THEORY AND THE PALEY–WIENER EXPANSION FOR FRACTIONAL BROWNIAN MOTION

By Kacha Dzhaparidze and Harry van Zanten

*Center for Mathematics and Computer Science and Vrije Universiteit Amsterdam*

In this paper we develop the spectral theory of the fractional Brownian motion (fBm) using the ideas of Krein's work on continuous analogues of orthogonal polynomials on the unit circle. We exhibit the functions which are orthogonal with respect to the spectral measure of the fBm and obtain an explicit reproducing kernel in the frequency domain. We use these results to derive an extension of the classical Paley–Wiener expansion of the ordinary Brownian motion to the fractional case.

**1. Introduction.** Let $X = (X_t)_{t \geq 0}$ be a fractional Brownian motion (fBm) with Hurst index $H \in (0, 1)$, that is, a continuous, centered Gaussian process with covariance function

$$\mathbb{E}X_s X_t = \tfrac{1}{2}(s^{2H} + t^{2H} - |s - t|^{2H}).$$

Say the process is defined on the probability space $(\Omega, \mathcal{F}, \mathbb{P})$, and for some fixed time horizon $T \geq 0$, define the linear space $\mathcal{H}_T$ as the closure in $L^2(\mathbb{P})$ of the (complex) linear span of a collection of random variables $\{X_t : t \in [0, T]\}$.

So-called linear problems for the fBm are problems in which it is required to find elements of the Hilbert space $\mathcal{H}_T$ with certain specific properties. In the 1960's the basic linear problems like prediction, interpolation, moving average representation etc. were treated by various authors; see Molchan (2003) for an overview of these contributions. However, these results did not become widely known. Many of them were rediscovered during the last decade when new application areas like telecommunication networks and mathematical finance stimulated a renewed interest in the fBm. Recent contributions dealing with linear problems can be found,

---

for instance, in Gripenberg and Norros (1996), Decreusefond and Üstünel (1999), Norros, Valkeila and Virtamo (1999), Nuzman and Poor (2000) and Pipiras and Taqqu (2001).

There exist several representations of the fBm that give insight into the structure of the linear space $\mathcal{H}_T$. An important example is the spectral representation

$$\mathbb{E} X_s X_t = \int_{\mathbb{R}} \frac{(e^{i\lambda t} - 1)(e^{-i\lambda s} - 1)}{\lambda^2} \mu(d\lambda) = \langle e_t, e_s \rangle_\mu,$$

where $\mu(d\lambda) = (2\pi)^{-1} \sin(\pi H) \Gamma(1 + 2H) |\lambda|^{1-2H} \, d\lambda$ is the spectral measure of the fBm and $e_t(\lambda) = (e^{i\lambda t} - 1)/i\lambda$ [see, e.g., Yaglom (1987), page 407 or Samorodnitsky and Taqqu (1994), page 328]. If we define $\mathcal{L}_T$ as the closure in $L^2(\mu)$ of the (complex) linear span of the collection of functions $\{e_t : t \in [0, T]\}$, this representation gives rise to an isometry between $\mathcal{H}_T$ and the function space $\mathcal{L}_T$, determined by the relation $X_t \longleftrightarrow e_t$. We can use this spectral isometry to reformulate a linear problem for the fBm in spectral terms. It then becomes a linear problem in the function space $\mathcal{L}_T$, which has the advantage that we have mathematical tools like Fourier-type techniques, at our disposal.

In this paper we present new results regarding the fine analytical structure of the frequency domain $\mathcal{L}_T$. In particular, we exhibit certain "orthogonal functions" with respect to the spectral measure $\mu$ of the fBm and we obtain an explicit reproducing kernel for $\mathcal{L}_T$, turning it into a reproducing kernel Hilbert space (RKHS). To illustrate the significance of these new frequency domain results for the fBm, we apply them to derive a generalization to the fractional case of the classical Paley–Wiener expansion of the ordinary Brownian motion [cf. Paley and Wiener (1934)]. In spectral terms, obtaining a series expansion translates to finding an orthonormal basis of the space $\mathcal{L}_T$. We achieve this by using the RKHS structure and the explicit expression that we have for the reproducing kernel.

It is well known that orthogonal polynomials on the unit circle are very useful in the spectral analysis of stationary time series. They can be used to solve problems like prediction and interpolation, and are also useful in connection with likelihood estimation and testing [see, e.g., Grenander and Szegö (1958)]. In a classical paper, Krein (1955) introduced certain continuous analogues of orthogonal polynomials on the unit circle. We refer to Akhiezer and Rybalko (1968) for a more elaborate treatment, including detailed proofs of Krein's statements. As was shown by Kailath, Vieira and Morf (1978), Krein's orthogonal functions play the same role in the spectral theory of continuous-time processes with stationary increments as the orthogonal polynomials do in times series theory. For a certain class of processes with stationary increments, Kailath, Vieira and Morf (1978) pointed out how the orthogonal functions describe the structure of the frequency domain. In the present paper we develop the spectral theory of the fBm along the same lines.

The results of Akhiezer and Rybalko (1968) and Kailath, Vieira and Morf (1978) highly depend on the "signal plus white noise" structure of the processes that they consider. It turns out, however, that Krein's ideas are also applicable for the fBm, which is not of the latter type. The key point is that the fBm can be "whitened" in the sense that integration of an appropriate deterministic kernel with respect to the fBm yields a continuous Gaussian martingale, the so-called "fundamental martingale." This was first proved in the 1960's by Molchan [cf. Molchan (2003)]. For alternative, more recent approaches see also Decreusefond and Üstünel (1999), Norros, Valkeila and Virtamo (1999), Nuzman and Poor (2000), Pipiras and Taqqu (2001) or Dzhaparidze and Ferreira (2002). Conversely, it is well known that by integrating a certain deterministic kernel with respect to the fundamental martingale, we can recover the fBm. Since these results play an important role in this paper, their precise statements will be recalled in Section 2.

The whitening and moving average formulas for the fBm provide us with the starting point for the development of the spectral theory. They give rise to a Hilbert space isometry $U$ between $\mathcal{L}_T$ and the space $L^2([0, T], V)$, where $V$ is the variance function of the fundamental martingale. In the ordinary Brownian case $H = 1/2$ this isometry is simply the Fourier transform. We will show that, for $H \neq 1/2$ it is also a Fourier-type integral transformation, and obtain an explicit expression for the Fourier kernel in terms of Bessel functions. Using this Fourier kernel, we then introduce a function $S_T$ that will turn out to be a reproducing kernel on $\mathcal{L}_T$. An explicit expression for this kernel will be derived in a number of steps. First we shall use the Bessel differential equation to prove that the properly normalized Fourier kernels satisfy Krein's continuous version of the recurrence relation for orthogonal polynomials. It will then be rather straightforward to obtain a Christoffel–Darboux-type formula for $S_T$. In combination with the expression for the Fourier kernel, this will lead to an explicit formula for the reproducing kernel on $\mathcal{L}_T$. This program is carried out in Sections 3–6.

In Sections 7 and 8 we use the new results on the structure of the space $\mathcal{L}_T$ to derive an extension to the fractional case of the classical Paley–Wiener expansion of the ordinary Brownian motion. We will first use the reproducing kernel to find a suitable orthogonal basis of $\mathcal{L}_T$. By transporting of this basis to the space $\mathcal{H}_T$, we will prove that the fBm admits the series expansion

$$\sum_{n \in \mathbb{Z}} \frac{e^{2i\omega_n t} - 1}{2i\omega_n} Z_n, \qquad t \in [0, 1],$$

where the $\omega_n$ are the real-valued zeros of the Bessel function $J_{1-H}$ and the $Z_n$ are independent, complex-valued Gaussian random variables with zero mean and a variance that can be expressed explicitly in terms of Bessel functions and their real zeros. Using the fact that $\sqrt{z}J_{1/2}(z) = \sqrt{2/\pi}\sin z$, it can be seen that for $H = 1/2$ this indeed reduces to

$$\sum_{n\in\mathbb{Z}}\frac{e^{2in\pi t}-1}{2in\pi}Z_n, \qquad t\in[0,1],$$

with the $Z_n$ i.i.d., standard Gaussian. This is the expression that Paley and Wiener (1934) used as the definition of the standard Brownian motion.

We will also briefly consider questions like the rate of convergence of the Paley–Wiener expansion, and possible extensions to the fractional Brownian sheet. In particular, we will argue that the expansion is rate-optimal, in the sense of Kühn and Linde (2002). This is obviously a desirable feature if the expansion is used for simulation purposes and is also relevant in connection with the small ball problem for the fractional Brownian sheet [see, e.g., Li and Linde (1999) and Li and Shao (2001) for the precise connections].

## 2. Auxiliary facts and notation.

The spectral representation can be used to define a stochastic integral with respect to $X$ of a large class of deterministic integrands. In this paper we denote the indicator function $\mathbb{1}_{(0,t)}$ of the interval $(0,t)$ simply by $\mathbb{1}_t$. Using this notation, as well as the previous notation $e_t(\lambda) = (\exp(i\lambda t)-1)/i\lambda$, we may write $e_t = \hat{\mathbb{1}}_t$. Here and elsewhere below we adopt the usual convention to denote the Fourier transform of a function $f \in L^2(\mathbb{R})$ by $\hat{f}$, that is,

$$\hat{f}(\lambda) = \int_{\mathbb{R}} f(x)e^{i\lambda x}\,dx.$$

Now consider the class of functions $\mathcal{I}_T = \{f \in L^2[0,T] : \hat{f} \in L^2(\mathbb{R}, \mathcal{B}(\mathbb{R}), \mu)\}$ and endow it with the inner product $\langle f, g\rangle_{\mathcal{I}_T} = \langle \hat{f}, \hat{g}\rangle_\mu$. Then the spectral representation can be written as $\mathbb{E}X_s X_t = \langle \mathbb{1}_s, \mathbb{1}_t\rangle_{\mathcal{I}_T}$. In particular, the mapping $\mathbb{1}_t \to X_t$ extends to a linear map $I : \mathcal{I}_T \to \mathcal{H}_T$ with the property that $I(\mathbb{1}_t) = X_t$ and for $f, g \in \mathcal{I}_T$,

$$\mathbb{E}I(f)\overline{I(g)} = \langle \hat{f}, \hat{g}\rangle_\mu.$$

We denote the random variable $I(f)$ by $\int f\,dX$ or $\int_0^T f(t)\,dX_t$, and call it the integral of $f$ with respect to $X$. We note that, in general, not every element of $\mathcal{H}_T$ can be represented as such an integral since for $H > 1/2$ the space $\mathcal{I}_T$ is not complete [see Pipiras and Taqqu (2001)].

Let us now introduce an integral with respect to $X$ which plays an important role in this paper. For $t \geq 0$, define the kernel $m_t$ by

$$(2.1) \quad m_t(u) = \frac{1}{2H\Gamma(H+1/2)\Gamma(3/2-H)}u^{1/2-H}(t-u)^{1/2-H}\mathbb{1}_t(u),$$

where $\Gamma$ denotes Euler's gamma function. Then, for every $t \in [0, T]$, it holds that $m_t \in \mathcal{I}_T$, and in view of Poisson's integral formula for the Bessel function [e.g., Watson (1944), Section 3.3], it is not hard to see that the Fourier transform $\hat{m}_t$ of $m_t$ is given by

$$(2.2) \qquad \hat{m}_t(\lambda) = \begin{cases} \dfrac{\sqrt{\pi}}{2H\Gamma(H+1/2)}\left(\dfrac{t}{\lambda}\right)^{1-H} e^{i\lambda t/2} J_{1-H}\left(\dfrac{\lambda t}{2}\right), & \lambda \neq 0, \\[2mm] \dfrac{\sqrt{\pi}}{2H\Gamma(H+1/2)\Gamma(2-H)2^{2-2H}} t^{2-2H}, & \lambda = 0, \end{cases}$$

where $J_{1-H}$ is the Bessel function of the first kind of order $1-H$ [for details see Dzhaparidze and Ferreira (2002), Proposition 2.2]. On evaluating $\hat{m}_t$ at $\lambda = 0$, one has to take into consideration the basic property

$$(2.3) \qquad z^{-\nu} J_\nu(z) \to \frac{1}{2^\nu \Gamma(\nu+1)} \qquad \text{as } z \to 0$$

of the Bessel function. For convenience, we introduce a special notation $d_H^2$ for the constant that occurs in the second line of (2.2). In the literature this constant is often given in an alternative form, namely,

$$(2.4) \qquad d_H^2 = \frac{\Gamma(3/2-H)}{2H\Gamma(H+1/2)\Gamma(3-2H)}.$$

The identity of the two expressions is a result of Legendre's duplication formula for the gamma function.

Now, for $t \in [0, T]$, we can consider the random variable

$$(2.5) \qquad M_t = \int m_t(u)\, dX_u = \int_0^t m_t(u)\, dX_u$$

in $\mathcal{H}_T$. As is proved in Dzhaparidze and Ferreira (2002), Theorem 2.3, it holds that $\mathbb{E}M_s M_t = \langle \hat{m}_s, \hat{m}_t \rangle_\mu = \hat{m}_{s \wedge t}(0)$. This shows that the process $M$ defined by (2.5) is a continuous Gaussian martingale with bracket $\langle M \rangle = \hat{m}(0)$. For convenience, this variance function will be denoted by $V$, so that due to (2.2), we have

$$(2.6) \qquad V_t = \mathbb{E}M_t^2 = d_H^2 t^{2-2H}$$

for all $t \geq 0$. Following Norros, Valkeila and Virtamo (1999), we call the process $M$ the fundamental martingale.

Next we recall the moving average representation of the fBm, which is the converse of (2.5). Let $x_t$ be defined by

$$x_t(u) = \left(t^{H-1/2}(t-u)^{H-1/2} - \int_u^t (t-v)^{H-1/2}\, dv^{H-1/2}\right) \mathbb{1}_t(u).$$

Then it holds that

$$(2.7) \qquad X_t = \int_0^t x_t(u)\, dM_u$$

for all $t \geq 0$, where $M$ is the fundamental martingale. More precisely, the process on the right-hand side defines an fBm with Hurst index $H$, so, in particular, we have that $x_t \in L^2([0, T], \mathcal{B}[0, T], V)$ for all $t \in [0, T]$ and

$$(2.8) \qquad \mathbb{E} X_s X_t = \int_0^{s \wedge t} x_s(u) x_t(u) \, dV_u = \langle x_s, x_t \rangle_V$$

for $s, t \in [0, T]$. We therefore consider the space $\mathcal{K}_T$, defined as the closure in $L^2([0, T], \mathcal{B}[0, T], V)$ of the (complex) linear span of the collection of functions $\{x_t : t \in [0, T]\}$. By construction, relation (2.8) shows that we have an isometry between $\mathcal{H}_T$ and $\mathcal{K}_T$, under which the correspondence $X_t \longleftrightarrow x_t$ holds true. Observe that under this isometry, we also have $M_t \longleftrightarrow \mathbb{1}_t$, so, for every $t \in [0, T]$, the indicator function $\mathbb{1}_t$ belongs to $\mathcal{K}_T$. It follows that $\mathcal{K}_T = L^2([0, T], \mathcal{B}[0, T], V)$.

In the remainder of the paper we write $L^2(\mu)$ and $L^2([0, T], V)$ instead of $L^2(\mathbb{R}, \mathcal{B}(\mathbb{R}), \mu)$ and $L^2([0, T], \mathcal{B}[0, T], V)$, respectively.

### 3. The transformation $\mathcal{K}_T \to \mathcal{L}_T$.

We have now associated three isometric Hilbert spaces with the fBm: the linear space $\mathcal{H}_T$, the frequency domain $\mathcal{L}_T$ and the space of integration kernels $\mathcal{K}_T$. The aforementioned isometries between $\mathcal{H}_T$ and $\mathcal{L}_T$ and between $\mathcal{H}_T$ and $\mathcal{K}_T$, determined by the relations $X_t \longleftrightarrow e_t$ and $X_t \longleftrightarrow x_t$, respectively, induce a direct isometry between the function spaces $\mathcal{K}_T$ and $\mathcal{L}_T$. We denote the map from $\mathcal{K}_T$ to $\mathcal{L}_T$ by $U$.

Our first result gives an explicit analytic description of the isometry $U : \mathcal{K}_T \to \mathcal{L}_T$. The theorem states that it is a Fourier-type integral transformation. The integration kernel is defined in terms of the function $\varphi : \mathbb{R} \to \mathbb{C}$, given by

$$(3.1) \quad \varphi(z) = \begin{cases} \Gamma(1-H)\left(\dfrac{z}{4}\right)^H e^{iz/2}\left(J_{-H}\left(\dfrac{z}{2}\right) + iJ_{1-H}\left(\dfrac{z}{2}\right)\right), & z \neq 0, \\ 1, & z = 0. \end{cases}$$

Here, as before, $\Gamma$ is Euler's gamma function and $J_\nu$ is the Bessel function of the first kind of order $\nu$. Observe that, in fact, $\varphi$ is defined on the whole complex plane. Moreover, property (2.3) of the Bessel function implies that $\varphi$ is an entire function. Evoking the well-known property

$$(3.2) \qquad \frac{d}{dz} z^\nu J_\nu(z) = z^\nu J_{\nu-1}(z)$$

of the Bessel function, one can easily see that the function $\varphi$ evaluated at $\lambda t$ and the earlier introduced Fourier transform $\hat{m}_t(\lambda)$ are related by the identity

$$(3.3) \qquad \varphi(\lambda t) = \frac{d\hat{m}_t(\lambda)}{d\hat{m}_t(0)}.$$

We need the following simple estimates for the function $\varphi$. The notation $a \lesssim b$ means that $a \leq cb$, where $c$ is positive constant that is universal or at least fixed throughout the paper.

LEMMA 3.1. *For every $\lambda \in \mathbb{R}$, the function $u \mapsto \varphi(u\lambda)$ belongs to $\mathcal{K}_T$ and its norm satisfies*

$$(3.4) \qquad \|\varphi(\cdot\lambda)\|_V \lesssim \begin{cases} 1 \wedge |\lambda|^{H-1/2}, & H \leq 1/2, \\ 1 \vee |\lambda|^{H-1/2}, & H > 1/2. \end{cases}$$

PROOF. The fact that $\varphi$ is analytic implies that it is bounded in a neighborhood of 0. Using also that $\sqrt{z} J_\nu(z)$ is bounded for $|z| \to \infty$, we see that for real $z$, $|\varphi(z)|$ is of order $|z|^{H-1/2}$ for large $|z|$. So for $H \leq 1/2$ we have $|\varphi(z)| \lesssim 1 \wedge |z|^{H-1/2}$, whence

$$\int_0^T |\varphi(u\lambda)|^2 \, dV_u \lesssim \int_0^T (u^{1-2H} \wedge |\lambda|^{2H-1}) \, du \lesssim 1 \wedge |\lambda|^{2H-1}.$$

For $H > 1/2$, it holds that $|\varphi(z)| \lesssim 1 \vee |z|^{H-1/2}$, so

$$\int_0^T |\varphi(u\lambda)|^2 \, dV_u \lesssim \int_0^T (u^{1-2H} \vee |\lambda|^{2H-1}) \, du \lesssim 1 \vee |\lambda|^{2H-1}.$$

This completes the proof. $\square$

The following theorem gives a complete description of the isometry $U : \mathcal{K}_T \to \mathcal{L}_T$.

THEOREM 3.2. *The linear transformation $U : \mathcal{K}_T \to \mathcal{L}_T$ is a Hilbert space isometry. For $f \in \mathcal{K}_T$, it holds that*

$$(3.5) \qquad Uf(\lambda) = \int_0^T f(u)\varphi(u\lambda) \, dV_u$$

*for $\mu$-almost all $\lambda \in \mathbb{R}$, where $\varphi$ is given by (3.1). The class of functions*

$$(3.6) \qquad \mathcal{L}_T' = \left\{ \psi \in \mathcal{L}_T : \int_{\mathbb{R}} \|\varphi(\cdot\lambda)\|_V |\psi(\lambda)| \mu(d\lambda) < \infty \right\}$$

*is dense in $\mathcal{L}_T$ and for $\psi \in \mathcal{L}_T'$ we have*

$$(3.7) \qquad U^*\psi(u) = U^{-1}\psi(u) = \int_{\mathbb{R}} \psi(\lambda)\overline{\varphi(u\lambda)}\mu(d\lambda)$$

*for $V$-almost every $u \in [0, T]$. Here $U^*$ denotes the adjoint of $U$.*

PROOF. By the Cauchy–Schwarz inequality and Lemma 3.1,

$$(3.8) \qquad \left| \int_0^T f(u) \varphi(u\lambda) \, dV_u \right| \leq \|\varphi(\cdot\lambda)\|_V \|f\|_V < \infty$$

for every $\lambda \in \mathbb{R}$. Hence, the right-hand side of (3.5) defines a linear transformation on $\mathcal{K}_T$. Let us denote this transformation by $A$.

Under the isometry $\mathcal{K}_T \to \mathcal{H}_T$, the indicator function $\mathbb{1}_t \in \mathcal{K}_T$ is mapped to the random variable $M_t \in \mathcal{H}_T$, given by (2.5). Under the spectral isometry, $M_t$ is mapped to the function $\lambda \mapsto \hat{m}_t(\lambda)$ in $\mathcal{L}_T$, given by (2.2). So to prove (3.5), we have to verify that the mapping $A$, defined on $\mathcal{K}_T$ by the right-hand side of (3.5), coincides with the isometry $U$ which is determined by the fact that $U\mathbb{1}_t = \hat{m}_t$ for $t \in [0, T]$.

By (3.3), we have

$$(3.9) \qquad \hat{m}_t(\lambda) = \int_0^t \varphi(u\lambda) \, dV_u.$$

So, indeed, $A\mathbb{1}_t = \hat{m}_t$ for every $t \in [0, T]$ and by linearity, $U$ coincides with $A$ on the set of simple functions in $\mathcal{K}_T$. Now take an arbitrary $f \in \mathcal{K}_T$. The simple functions are dense in $\mathcal{K}_T$, whence we can choose a sequence $f_n$ of simple functions such that $f_n \to f$ in $\mathcal{K}_T$. Then since $U$ is an isometry, we have $Af_n = Uf_n \to Uf$ in $\mathcal{L}_T \subseteq L^2(\mu)$. On the other hand, (3.8) implies that $Af_n \to Af$ pointwise on $\mathbb{R}$. But then $Af$ and $Uf$ must coincide for $\mu$-almost all $\lambda \in \mathbb{R}$, which proves (3.5).

Both the isometry $\mathcal{K}_T \to \mathcal{H}_T$ and the spectral isometry $\mathcal{H}_T \to \mathcal{L}_T$ preserve inner products, so the same holds for their composition $U \colon \mathcal{K}_T \to \mathcal{L}_T$. This implies that $U$ is unitary, that is, that $U^{-1} = U^*$, where $U^* \colon \mathcal{L}_T \to \mathcal{K}_T$ is the adjoint of $U$, determined by the relation $\langle Uf, \psi \rangle_\mu = \langle f, U^*\psi \rangle_V$ for all $f \in \mathcal{K}_T$ and $\psi \in \mathcal{L}_T$. Using (3.5), we see that, for $\psi \in \mathcal{L}_T$, we have

$$\langle Uf, \psi \rangle_\mu = \int_\mathbb{R} Uf(\lambda) \overline{\psi(\lambda)} \mu(d\lambda)$$

$$= \int_\mathbb{R} \left( \int_0^T f(u) \varphi(u\lambda) \, dV_u \right) \overline{\psi(\lambda)} \mu(d\lambda).$$

For $\psi \in \mathcal{L}'_T$, we may interchange the integrals, since by the preceding lemma and Cauchy–Schwarz,

$$\int_\mathbb{R} \left( \int_0^T |f(u)| |\varphi(u\lambda)| \, dV_u \right) |\psi(\lambda)| \mu(d\lambda) \leq \|f\|_V \int_\mathbb{R} \|\varphi(\cdot\lambda)\|_V |\psi(\lambda)| \mu(d\lambda),$$

which is finite by definition of $\mathcal{L}'_T$. It follows that

$$\langle Uf, \psi \rangle_\mu = \int_0^T f(u) \left( \int_\mathbb{R} \overline{\psi(\lambda)} \varphi(u\lambda) \mu(d\lambda) \right) dV_u,$$

which proves (3.7).

It remains to prove that $\mathcal{L}_T'$ is dense in $\mathcal{L}_T$. Let the $\mathcal{S}$ be the Schwarz space of rapidly decreasing functions on $\mathbb{R}$, that is, $C^\infty$-functions $f$ on $\mathbb{R}$ such that for all $m, n$, the derivative $f^{(n)}$ satisfies $|x|^m |f^{(n)}(x)| \to 0$ as $|x| \to \infty$, and let $\mathcal{S}_T$ be the space of Schwarz functions with support in $[0, T]$. By the preceding lemma, it clearly holds that $\mathcal{S} \cap \mathcal{L}_T \subseteq \mathcal{L}_T'$, so it is enough to show that $\mathcal{S} \cap \mathcal{L}_T$ is dense in $\mathcal{L}_T$.

Fix $t \in [0, T]$ and choose a sequence $w_n$ of $C^\infty$ probability densities such that $\mathrm{supp}(w_n) \subseteq [0, t]$ and such that the associated probability measures converges weakly to the dirac measure $\delta_t$ concentrated at $t$. Define

$$f_n(u) = \int_0^t \mathbb{1}_s(u) w_n(s) \, ds = \begin{cases} \int_u^t w_n(s) \, ds, & u \leq t, \\ 0, & u \in (t, T]. \end{cases}$$

Then $f_n$ is a $C^\infty$-function with compact support, so $f_n \in \mathcal{S}_T$. For fixed $\lambda \in \mathbb{R}$, the function $s \mapsto \hat{\mathbb{1}}_s(\lambda)$ is bounded and continuous, whence the weak convergence implies that

$$\hat{f}_n(\lambda) = \int_0^t \hat{\mathbb{1}}_s(\lambda) w_n(s) \, ds \to \int_0^t \hat{\mathbb{1}}_s(\lambda) \delta_t(ds) = \hat{\mathbb{1}}_t(\lambda)$$

for every $\lambda \in \mathbb{R}$. Observe also that, for $\lambda \in \mathbb{R}$,

$$|\hat{f}_n(\lambda)|^2 \leq \int_0^t |\hat{\mathbb{1}}_s(\lambda)|^2 w_n(s) \, ds \lesssim 1 \wedge \frac{1}{\lambda^2}.$$

By dominated convergence, it follows that $\hat{f}_n \to \hat{\mathbb{1}}_t$ in $L^2(\mu)$. Since the Fourier transform maps $\mathcal{S}_T$ into $\mathcal{S}$, the functions $\hat{f}_n$ belong to $\mathcal{S} \cap \mathcal{L}_T$. So for every $t \in [0, T]$, $\hat{\mathbb{1}}_t$ is the $L^2(\mu)$-limit of a sequence of functions in $\mathcal{S} \cap \mathcal{L}_T$. Since $\mathcal{L}_T$ is the closure in $L^2(\mu)$ of the linear span of the functions $\hat{\mathbb{1}}_t$, $t \in [0, T]$, this shows that $\mathcal{S} \cap \mathcal{L}_T$ is, indeed, dense in $\mathcal{L}_T$. $\square$

Observe that since we have

$$(3.10) \qquad J_{1/2}(z) = \sqrt{\frac{2}{\pi z}} \sin z, \qquad J_{-1/2}(z) = \sqrt{\frac{2}{\pi z}} \cos z, \qquad z \neq 0$$

and $\Gamma(1/2) = \sqrt{\pi}$, it holds that $\varphi(z) = e^{iz}$ in the standard Brownian motion case $H = 1/2$. So in that case, the map $U$ is simply the Fourier transform. For general $H \in (0, 1)$, we can view it as a fractional version of the Fourier transform.

It seems worth mentioning that fractional integration theory enters in the present context via the simple observation that the Fourier transform $\hat{m}_t$ of the kernel $m_t$, defined by (2.1), is expressible in terms of the fractional integral of order $3/2 - H$ of the function $u \mapsto u^{1/2-H} \exp(iu\lambda)$. This can be seen by comparing (2.2) with formula 9.1.10 of Samko, Kilbas and Marichev (1993). Specifically, we have that

$$\hat{m}_t(\lambda) = \frac{1}{\Gamma(H + 1/2)} I_{0+}^{3/2-H}(u^{1/2-H} e^{iu\lambda})(t).$$

By (3.3), it follows that, for the Fourier-kernel of the map $U$, it holds that

$$V_t' \varphi(\lambda t) = \frac{1}{\Gamma(H + 1/2)} I_{0+}^{1/2-H}(u^{1/2-H} e^{iu\lambda})(t).$$

Hence, using fractional integration by parts, we see that, for $f \in \mathcal{K}_T$,

$$Uf(\lambda) = \frac{1}{\Gamma(H + 1/2)} \int_0^T t^{1/2-H} e^{i\lambda t} I_{T-}^{1/2-H} f(t) \, dt$$

$$= \frac{1}{\Gamma(H + 1/2)} \mathcal{F}(u^{1/2-H} I_{T-}^{1/2-H} f(u))(\lambda),$$

provided, of course, that the fractional integral of order $1/2 - H$ of $f$ and the Fourier transform (denoted by $\mathcal{F}$) exist. The composition rule of fractional integration operators implies that, for $\psi \in \mathcal{L}_T$,

$$U^{-1} \psi(t) = \Gamma(H + 1/2) I_{T-}^{H-1/2}(u^{H-1/2} \mathcal{F}^{-1} \psi(u))(t).$$

Note, for instance, that for a deterministic integrand $f \in \mathcal{I}_T$, the latter expression for $U^{-1}$, in combination with the spectral isometry, yields the relation

$$\int_0^T f(u) \, dX_u = \Gamma(H + 1/2) \int_0^T I_{T-}^{H-1/2}(u^{H-1/2} f(u))(t) \, dM_t,$$

where $M$ is the fundamental martingale. For $f = \mathbb{1}_t$, this reduces to the moving average representation (2.7). In general, the expression of the operators $U$ and $U^{-1}$, in terms of Riemann–Liouville operators, can be very useful for the evaluation of the transforms in concrete cases, since many explicit formulas for fractional integrals are known. We will, however, not need this connection in the present paper. The proofs of our results do not use any fractional calculus.

Relation (3.5) gives an analytic description of the functions in $\mathcal{L}_T$. In particular, it allows us to prove that every function in $\mathcal{L}_T$ is the restriction to $\mathbb{R}$ of an entire function. Strictly speaking, the elements of $\mathcal{L}_T$ are, of course, equivalence classes of functions. Two functions represent the same element if they coincide $\mu$-almost everywhere. Theorem 3.2 implies that every equivalence class can be represented by an entire function.

COROLLARY 3.3. *Every element of $\mathcal{L}_T$ has a version that is the restriction to $\mathbb{R}$ of an entire function.*

Proof. For $f \in \mathcal{K}_T$, consider the complex function

$$(3.11) \qquad z \mapsto \int_0^T f(u)\varphi(uz)\,dV_u.$$

Since $\varphi$ is entire, this function is well defined and easily seen to be continuous on $\mathbb{C}$. To prove that the function is analytic, consider a closed path $\gamma$ in the complex plane. By Fubini's theorem and Cauchy's theorem,

$$\oint_\gamma \left( \int_0^T f(u)\varphi(uz)\,dV_u \right) dz = \int_0^T f(u) \left( \oint_\gamma \varphi(uz)\,dz \right) dV_u = 0.$$

Hence, by Morera's theorem, the function defined by (3.11) is entire. □

In the remainder of the paper, if we consider an element $\psi \in \mathcal{L}_T$, we will always assume this to be the smooth version.

## 4. The reproducing kernel on $\mathcal{L}_T$.

As was established in Theorem 3.2, the transformation $U$ is of Fourier-type, generated by the Fourier kernel $\varphi$. This motivates us to introduce the function $S_T$ on $\mathbb{R} \times \mathbb{R}$, defined by

$$S_T(\omega, \lambda) = \int_0^T \overline{\varphi(u\omega)}\varphi(u\lambda)\,dV_u,$$

where $\varphi$ is, of course, given by (3.1) again, and $V$ by (2.6). The Cauchy–Schwarz inequality and Lemma 3.1 imply that $S_T$ is well defined. Moreover, (3.5) implies that, for fixed $\omega \in \mathbb{R}$, the function $\lambda \mapsto S_T(\omega, \lambda)$ is the image under $U$ of the function $u \mapsto \overline{\varphi(u\omega)}$. In particular, we see that $\lambda \mapsto S_T(\omega, \lambda)$ belongs to $\mathcal{L}_T$ for every $\omega \in \mathbb{R}$. It clearly holds that $S_T(\omega, \lambda) = \overline{S_T(\lambda, \omega)}$. Since $\varphi(0) = 1$, relation (3.9) implies that $S_T(0, 0) = V_T$ and $S_T(0, \lambda) = \hat{m}_T(\lambda)$.

The following theorem states that $S_T$ acts as a reproducing kernel on the spectral space $\mathcal{L}_T$, turning it into an RKHS.

THEOREM 4.1. For every $\psi \in \mathcal{L}_T$, we have

$$\int_\mathbb{R} \psi(\lambda)\overline{S_T(\omega, \lambda)}\mu(d\lambda) = \psi(\omega),$$

for all $\omega \in \mathbb{R}$.

Proof. Suppose first that $\psi \in \mathcal{S} \cap \mathcal{L}_T$, where $\mathcal{S}$ is the Schwarz space of rapidly decreasing functions. Then by Fubini's theorem and Theorem 3.2,

$$\int_\mathbb{R} \psi(\lambda)\overline{S_T(\omega, \lambda)}\mu(d\lambda) = \int_\mathbb{R} \psi(\lambda)\left( \int_0^T \varphi(u\omega)\overline{\varphi(u\lambda)}\,dV_u \right)\mu(d\lambda)$$

$$= \int_0^T \varphi(u\omega)\left( \int_\mathbb{R} \psi(\lambda)\overline{\varphi(u\lambda)}\mu(d\lambda) \right) dV_u$$

$$= \int_0^T \varphi(u\omega) U^{-1} \psi(u) \, dV_u$$

$$= U(U^{-1}\psi(\omega)) = \psi(\omega)$$

for $\mu$-almost all $\omega \in \mathbb{R}$. The interchanging of the integration order is justified by the fact that $\psi$ is rapidly decreasing.

Now let $\psi \in \mathcal{L}_T$ be arbitrary. Since $\mathcal{S} \cap \mathcal{L}_T$ is dense in $\mathcal{L}_T$ (see the proof of Theorem 3.2), we can choose functions $\psi_n \in \mathcal{S} \cap \mathcal{L}_T$ such that $\psi_n \to \psi$ in $L^2(\mu)$. By the remarks preceding the theorem, the function $\lambda \mapsto S_T(\omega, \lambda)$ belongs to $L^2(\mu)$ for fixed $\omega \in \mathbb{R}$. So by Cauchy–Schwarz, we have

$$\left| \int_{\mathbb{R}} \psi_n(\lambda) \overline{S_T(\omega, \lambda)} \mu(d\lambda) - \int_{\mathbb{R}} \psi(\lambda) \overline{S_T(\omega, \lambda)} \mu(d\lambda) \right|$$

$$\leq \|\psi_n - \psi\|_\mu \|S_T(\omega, \cdot)\|_\mu \to 0$$

for every $\omega \in \mathbb{R}$. By the preceding paragraph, the first integral on the left-hand side equals $\psi_n(\omega)$ for $\mu$-almost every $\omega$. So the functions $\psi_n$ converge $\mu$-almost everywhere to the function

$$\omega \mapsto \int_{\mathbb{R}} \psi(\lambda) \overline{S_T(\omega, \lambda)} \mu(d\lambda),$$

and they converge in $L^2(\mu)$ to $\psi$. Then the two limits must coincide $\mu$-almost everywhere. Since both functions are continuous (Corollary 3.3), the proof is complete. □

A first simple consequence of the RKHS structure is that the "kernel functions" $\lambda \mapsto S_T(\omega, \lambda)$ span $\mathcal{L}_T$.

COROLLARY 4.2.    *The space $\mathcal{L}_T$ is the closure in $L^2(\mu)$ of the linear span of the collection of functions $\{\lambda \mapsto S_T(\omega, \lambda) : \omega \in \mathbb{R}\}$.*

PROOF.    We already noted in the beginning of the section that every function $S_T(\omega, \cdot)$ belongs to $\mathcal{L}_T$, so the closure of the linear span of $\{\lambda \mapsto S_T(\omega, \lambda) : \omega \in \mathbb{R}\}$ is certainly contained in $\mathcal{L}_T$.

To prove the inclusion in the other direction, take $\psi \in \mathcal{L}_T$ and suppose that $\psi$ is orthogonal to every kernel function $S_T(\omega, \cdot)$, that is,

$$\int_{\mathbb{R}} \psi(\lambda) \overline{S_T(\omega, \lambda)} \mu(d\lambda) = 0$$

for all $\omega \in \mathbb{R}$. Then, by the reproducing property of $S_T$, we see that $\psi$ vanishes $\mu$-almost everywhere. □

It seems useful to briefly discuss the relation between the preceding frequency-domain results and the so-called "time-domain RKHS." The latter space is constructed by associating to every element $H \in \mathcal{H}_T$ a function $t \mapsto \mathbb{E} H X_t$ on $[0, T]$. These functions are the elements of the time-domain RKHS and the inner product of two functions $t \mapsto \mathbb{E} H X_t$ and $t \mapsto \mathbb{E} H' X_t$ is defined as $\mathbb{E} H H'$. By construction, the resulting Hilbert space is isometric to $\mathcal{H}_T$, and the covariance function $r(s, t) = \mathbb{E} X_s X_t$ is the reproducing kernel on the space.

The following theorem clarifies the relation between the two reproducing kernel Hilbert space structures.

THEOREM 4.3. *The time-domain RKHS is the closure of the linear span of the collection of functions $\{t \mapsto \hat{\mathbb{1}}_t(\lambda) : \lambda \in \mathbb{R}\}$ with respect to the inner product*

$$\langle t \mapsto \hat{\mathbb{1}}_t(\omega), t \mapsto \hat{\mathbb{1}}_t(\lambda) \rangle = S_T(\omega, \lambda).$$

PROOF. The spectral isometry shows that the time-domain RKHS is given by all functions $t \mapsto \langle \psi, \hat{\mathbb{1}}_t \rangle_\mu$ on $[0, T]$, where $\psi$ runs through $\mathcal{L}_T$, the inner product of two elements $t \mapsto \langle \psi, \hat{\mathbb{1}}_t \rangle_\mu$ and $t \mapsto \langle \xi, \hat{\mathbb{1}}_t \rangle_\mu$ being given by $\langle \psi, \xi \rangle_\mu$. By Corollary 4.2, it follows that the time-domain RKHS is the closure of the linear span of collection of functions $\{t \mapsto \langle S_T(\lambda, \cdot), \hat{\mathbb{1}}_t \rangle_\mu : \lambda \in \mathbb{R}\}$. By the reproducing property, we have $\langle S_T(\lambda, \cdot), \hat{\mathbb{1}}_t \rangle_\mu = \hat{\mathbb{1}}_t(\lambda)$ and the inner product between the functions $t \mapsto \hat{\mathbb{1}}_t(\lambda)$ and $t \mapsto \hat{\mathbb{1}}_t(\omega)$ equals $S_T(\omega, \lambda)$. □

We note that the moving average representation of the fBm implies that the inner products in the time-domain RKHS can be expressed in terms of Riemann–Liouville fractional integration operators [cf., e.g., Hult (2003)]. Theorem 4.3 thus yields an expression for the reproducing kernel $S_T$ of $\mathcal{L}_T$ in terms of fractional integrals. In the present paper this connection plays no further role, however, and we will now return to the frequency domain analysis.

**5. Differential equations.** Our next intension is to obtain an explicit analytic expression for the reproducing kernel $S_T$. This goal is achieved in Corollary 6.2 below. In this section we will show that the Fourier kernel $\varphi$ is subject to a "recursion relationship" that is also encountered in the study of other types of processes with stationary increments [see, e.g., Kailath, Vieira and Morf (1978), where a short account can be found of the classical result of Krein (1955) concerning the "signal plus white noise" model].

The Bessel function $J_\nu$ satisfies the second order ordinary differential equation

$$(5.1) \qquad J_\nu''(z) + \frac{1}{z} J_\nu'(z) + \left(1 - \frac{\nu^2}{z^2}\right) J_\nu(z) = 0.$$

In view of the representation (2.2), this yields the following differential equation for the functions $\hat{m}_t \in \mathcal{L}_T$.

LEMMA 5.1.  *For every $\lambda \in \mathbb{R}$, we have*

$$\frac{\partial^2 \hat{m}_t(\lambda)}{\partial t^2} = i\lambda \frac{\partial \hat{m}_t(\lambda)}{\partial t} - 2\gamma_t \left( \frac{\partial \hat{m}_t(\lambda)}{\partial t} - \frac{i\lambda}{2} \hat{m}_t(\lambda) \right),$$

*where $\gamma_t = (H - 1/2)/t$.*

PROOF.   By (2.2), we have $\hat{m}_t(\lambda) = c_\lambda f_\nu(z) J_\nu(z)$, where $c_\lambda$ is a constant (depending on $\lambda$), $z = \lambda t/2$, $\nu = 1 - H$ and $f_\nu(z) = z^\nu \exp(iz)$. It is easily verified that

$$\frac{\partial}{\partial t} f_\nu(z) = \frac{f_\nu(z)(\nu + iz)}{t}, \qquad \frac{\partial^2}{\partial t^2} f_\nu(z) = \frac{f_\nu(z)((\nu + iz)^2 - \nu)}{t^2}.$$

The differential equation is now a straightforward consequence of (5.1).   □

This differential equation for the functions $t \mapsto \hat{m}_t(\lambda)$ gives rise to a differential equation for the Fourier kernel $\varphi(t\lambda)$. We present this as a system of equations for the function $P(t, \lambda)$, defined for $t > 0$ by

$$(5.2) \qquad P(t, \lambda) = \varphi(t\lambda) \sqrt{dV_t/dt},$$

and its reciprocal

$$(5.3) \qquad P^*(t, \lambda) = e^{i\lambda t} \overline{P(t, \lambda)}.$$

Note that

$$(5.4) \qquad S_T(\omega, \lambda) = \int_0^T \overline{P(t, \lambda)} P(t, \lambda) \, dt.$$

We also observe that since the functions $\hat{m}_t$ correspond to the fundamental martingale $M$ under the spectral isometry, the functions

$$\int_0^t P(s, \cdot) \, ds$$

correspond to an ordinary Brownian motion. In particular, we have that

$$\int_{\mathbb{R}} \int_0^s P(u, \lambda) \, du \int_0^t \overline{P(v, \lambda)} \, dv \, \mu(d\lambda) = s \wedge t.$$

If we interchange the integrals on the left-hand side and differentiate, we get the formal expression

$$\int_{\mathbb{R}} P(s,\lambda)\overline{P(t,\lambda)}\mu(d\lambda) = \delta(t-s).$$

In this sense, the functions $P(t,\cdot)$ are orthogonal with respect to the spectral measure $\mu$.

The following theorem shows that the orthogonal functions satisfy Krein's continuous analogue of the usual recurrence formulas for orthogonal polynomials on the unit circle.

THEOREM 5.2. *For every $\lambda \in \mathbb{R}$, the function $P(t,\lambda)$ and its reciprocal $P^*(t,\lambda)$, defined by (5.2) and (5.3), satisfy the equations*

$$\frac{\partial P(t,\lambda)}{\partial t} = i\lambda P(t,\lambda) - \gamma_t P^*(t,\lambda)$$

*and*

$$\frac{\partial P^*(t,\lambda)}{\partial t} = -\gamma_t P(t,\lambda),$$

*where $\gamma_t = (H-1/2)/t$.*

PROOF. First of all, let us express $P(t,\cdot)$ and $P^*(t,\cdot)$ in terms of the function $\hat{m}_t$ given by (2.2). We have $d\hat{m}_t(\lambda) = \varphi(t\lambda)\,dV_t$ and $V_t = d_H^2 t^{2-2H}$ [see (3.3) and (2.6)], hence,

$$(5.5) \qquad P(t,\lambda) = \frac{t^{H-1/2}\hat{m}_t'(\lambda)}{\sqrt{(2-2H)d_H^2}},$$

where the prime denotes differentiation with respect to $t$. It follows that

$$(5.6) \qquad P^*(t,\lambda) = \frac{t^{H-1/2}(\hat{m}_t')^*(\lambda)}{\sqrt{(2-2H)d_H^2}},$$

where $(\hat{m}_t')^*(\lambda) = \exp(i\lambda t)\overline{\hat{m}_t'(\lambda)}$ is the reciprocal of $\hat{m}_t'$. Observe that

$$\hat{m}_t'(\lambda) - \frac{i\lambda}{2}\hat{m}_t(\lambda) = e^{i\lambda t/2}\frac{d}{dt}(e^{-i\lambda t/2}\hat{m}_t(\lambda)),$$

so that

$$\left(\hat{m}_t'(\lambda) - \frac{i\lambda}{2}\hat{m}_t\right)^*(\lambda) = e^{i\lambda t}\overline{e^{i\lambda t/2}\frac{d}{dt}(e^{-i\lambda t/2}\hat{m}_t(\lambda))}$$

$$= e^{i\lambda t/2}\frac{d}{dt}(e^{-i\lambda t/2}\hat{m}_t(\lambda))$$

$$= \hat{m}_t'(\lambda) - \frac{i\lambda}{2}\hat{m}_t(\lambda).$$

Since $\hat{m}_t$ is self-reciprocal, that is, $\hat{m}_t^*(\lambda) = \exp(i\lambda t)\overline{\hat{m}_t(\lambda)} = \hat{m}_t(\lambda)$, the latter identity implies that $(\hat{m}_t')^*(\lambda) = \hat{m}_t'(\lambda) - i\lambda\hat{m}_t(\lambda)$. Combining this with (5.6), we find that

$$(5.7) \qquad P^*(t,\lambda) = \frac{t^{H-1/2}(\hat{m}_t'(\lambda) - i\lambda\hat{m}_t(\lambda))}{\sqrt{(2-2H)d_H^2}}.$$

The first statement of the theorem now follows from differentiation of (5.5), taking Lemma 5.1 and (5.5) and (5.7) into account. Similarly, the second statement is obtained by differentiating (5.7). $\square$

**6. Christoffel–Darboux formula.** As in the theory of orthogonal polynomials and their continuous analogues, the "recurrence relations" presented in Theorem 5.2 allow us to derive a closed-form expression for the reproducing kernel $S_T$ [cf., e.g., Grenander and Szegö (1958), Section 2.3, and Kailath, Vieira and Morf (1978), formula (48)].

THEOREM 6.1. *Let $P(t,\lambda)$ and its reciprocal $P^*(t,\lambda)$ be defined by (5.2) and (5.3). For all $T > 0$ and $\omega, \lambda \in \mathbb{R}$, we have*

$$(6.1) \qquad i(\lambda - \omega)S_T(\omega,\lambda) = \overline{P(T,\omega)}P(T,\lambda) - \overline{P^*(T,\omega)}P^*(T,\lambda).$$

PROOF. We view the left-hand side and right-hand side of (6.1) as functions in $T$. Recall that we have (5.4). Using Theorem 5.2, a straightforward calculation shows that

$$i(\lambda - \omega)\overline{P(t,\omega)}P(t,\lambda) = \frac{\partial}{\partial t}(\overline{P(t,\omega)}P(t,\lambda) - \overline{P^*(t,\omega)}P^*(t,\lambda)).$$

This shows that the functions on both sides of (6.1) have the same derivative with respect to $T$, which implies that their difference is independent of $T$. So for every $T > 0$, we have

$$(6.2) \qquad i(\lambda - \omega)S_T(\omega,\lambda) = \overline{P(T,\omega)}P(T,\lambda) - \overline{P^*(T,\omega)}P^*(T,\lambda) + C(\omega,\lambda)$$

for some constant $C(\omega,\lambda)$, and it remains to show that $C(\omega,\lambda) = 0$.

For $H < 1/2$, the functions $t \mapsto P(t,\lambda)$ are bounded for every $\lambda \in \mathbb{R}$ and it holds that $P(0,\lambda) = 0$. So in this case we can let $T \to 0$ in (6.2) to see that $C(\omega,\lambda) = 0$. Since the Bessel function $J_\nu$ is analytic in $\nu$ [see Watson (1944), page 44], (6.2) shows that as a function in $H$, $C(\omega,\lambda)$ can be extended to an analytic function on the open disc of diameter 1, centered at $1/2$. We just saw that it vanishes for all $H$ in the interval $(0, 1/2)$, whence a standard result from complex function theory implies that it vanishes on the entire disc [cf. Rudin (1987), Theorem 10.19]. In particular, $C(\omega,\lambda) = 0$ for all $H \in (0,1)$. $\square$

In combination with the explicit expression that we have for the orthogonal functions $\varphi(t\lambda)$, the preceding theorem yields an explicit analytic expression for the reproducing kernel $S_T$.

COROLLARY 6.2.  *The reproducing kernel $S_T$ admits the following representation*:

(i)  *For $\omega \neq \lambda$,*

$$\frac{S_T(2\omega, 2\lambda)}{S_T(0,0)} = (2 - 2H)\Gamma^2(1 - H)\left(\frac{T^2\omega\lambda}{4}\right)^H e^{iT(\lambda-\omega)}$$

$$\times \frac{J_{-H}(T\omega)J_{1-H}(T\lambda) - J_{1-H}(T\omega)J_{-H}(T\lambda)}{T(\lambda - \omega)}.$$

(ii)  *For $\omega \in \mathbb{R}$,*

$$\frac{S_T(2\omega, 2\omega)}{S_T(0,0)} = (2 - 2H)\Gamma^2(1 - H)\left(\frac{T\omega}{2}\right)^{2H}$$

$$\times \left(J_{1-H}^2(T\omega) + \frac{2H-1}{T\omega}J_{-H}(T\omega)J_{1-H}(T\omega) + J_{-H}^2(T\omega)\right).$$

PROOF.  Part (i) follows by straightforward calculations from the preceding theorem, the definition (5.2) of $P(t,\lambda)$, the explicit expression (3.1) for $\varphi$ and the fact that $dV_t/dt = (2 - 2H)V_t/t$.

To prove part (ii), we note that $S_t$ is analytic and, in particular, continuous, so we may derive an expression for $S_t(\omega,\omega)$ by letting $\lambda \to \omega$ in the expression that we found in part (i). It suffices to observe that as $\lambda \to \omega$, we have

$$\frac{J_{-H}(t\omega)J_{1-H}(t\lambda) - J_{1-H}(t\omega)J_{-H}(t\lambda)}{\lambda - \omega}$$

$$= J_{-H}(t\omega)\frac{J_{1-H}(t\lambda) - J_{1-H}(t\omega)}{\lambda - \omega} - J_{1-H}(t\omega)\frac{J_{-H}(t\lambda) - J_{-H}(t\omega)}{\lambda - \omega}$$

$$\to J_{-H}(t\omega)\frac{\partial}{\partial\omega}J_{1-H}(t\omega) - J_{1-H}(t\omega)\frac{\partial}{\partial\omega}J_{-H}(t\omega)$$

$$= t\left(J_{1-H}^2(t\omega) + \frac{2H-1}{t\omega}J_{-H}(t\omega)J_{1-H}(t\omega) + J_{-H}^2(t\omega)\right).$$

In the last step we have used the recurrence formulae

(6.3)
$$\frac{d}{dz}J_\nu(z) = \frac{\nu}{z}J_\nu(z) - J_{\nu+1}(z)$$

and

(6.4)
$$J_{\nu+2}(z) = \frac{2\nu + 2}{z}J_{\nu+1}(z) - J_\nu(z)$$

[see Watson (1944), page 45].  □

**7. Orthonormal basis in $\mathcal{L}_T$.** By the reproducing property, the inner product of the kernel functions $\lambda \mapsto S_T(2\omega, \lambda)$ and $\lambda \mapsto S_T(2\omega', \lambda)$ in $\mathcal{L}_T$ is given by $S_T(2\omega, 2\omega')$. Hence, by Corollary 6.2, these functions are orthogonal in $\mathcal{L}_T$ if $T\omega$ and $T\omega'$ are different zeros of $J_{1-H}$. In this section we prove that if we let $T\omega$ range over all zeros of $J_{1-H}$, we obtain an orthogonal basis of $\mathcal{L}_T$.

We first recall some facts about the zeros of the Bessel function of the first kind [see, e.g., Erdélyi, Magnus, Oberhettinger and Tricomi (1953), Section 7.9]. For $\nu > -1$, the Bessel function $J_\nu$ has a countable number of positive zeros that can be ordered according to magnitude. We denote them by $\lambda_{\nu,1} < \lambda_{\nu,2} < \cdots$. For positive $\nu$, the function $J_\nu$ satisfies $J_\nu(0) = 0$ and its negative zeros are given by $-\lambda_{\nu,1} > -\lambda_{\nu,2} > \cdots$. Hence, for $\nu \geq 0$, the zeros of $J_\nu$ can be ordered as $\cdots < \lambda_{\nu,-1} < \lambda_{\nu,0} = 0 < \lambda_{\nu,1} < \cdots$. To prove the completeness of the system of orthogonal functions, we need the following lemma, which is a consequence of the multiplicative decomposition of the Bessel function, or, more precisely, of the formula

$$(7.1) \qquad \frac{J_{\nu+1}(z)}{J_\nu(z)} = \sum_{k=1}^{\infty} \frac{2z}{\lambda_{\nu,k}^2 - z^2};$$

see Erdélyi, Magnus, Oberhettinger and Tricomi [(1953), formula 7.9.3].

LEMMA 7.1. *Let $J_\nu$ be a Bessel function of the first kind of nonnegative order $\nu \geq 0$ and let $\cdots < \lambda_{-1} < \lambda_0 = 0 < \lambda_1 < \cdots$ be its real zeros. Then*

$$\frac{J_{\nu+1}(z)}{J_\nu(z)} - \frac{J_{\nu+1}(w)}{J_\nu(w)} = \sum_{n \neq 0} \frac{z-w}{(\lambda_n - z)(\lambda_n - w)}$$

*for all $z, w \in \mathbb{R}$, for which the expressions are well defined.*

PROOF. Since the zeros of $J_\nu$ satisfy $\lambda_{-n} = -\lambda_n$ for all $n \in \mathbb{N}$, it holds that

$$\sum_{n \neq 0} \frac{z-w}{(\lambda_n - z)(\lambda_n - w)} = \sum_{n \in \mathbb{N}} \left( \frac{z-w}{(\lambda_n - z)(\lambda_n - w)} + \frac{z-w}{(\lambda_n + z)(\lambda_n + w)} \right).$$

Now observe that we have the identity

$$\frac{a-b}{(c-a)(c-b)} + \frac{a-b}{(c+a)(c+b)} = \frac{2a}{c^2 - a^2} - \frac{2b}{c^2 - b^2},$$

for all $a, b, c \in \mathbb{R}$ for which the expression makes sense. It follows that

$$\sum_{n \neq 0} \frac{z-w}{(\lambda_n - z)(\lambda_n - w)} = \sum_{n \in \mathbb{N}} \frac{2z}{\lambda_n^2 - z^2} - \sum_{n \in \mathbb{N}} \frac{2w}{\lambda_n^2 - w^2}.$$

By (7.1), the right-hand side is equal to

$$\frac{J_{\nu+1}(z)}{J_\nu(z)} - \frac{J_{\nu+1}(w)}{J_\nu(w)},$$

and the proof of the lemma is complete. $\square$

We can now present the orthonormal basis in $\mathcal{L}_T$ and the associated expansion formula.

THEOREM 7.2. *Let* $\cdots < \omega_{-1} < \omega_0 = 0 < \omega_1 < \cdots$ *be the real zeros of* $J_{1-H}$ *and, for* $n \in \mathbb{Z}$, *define the function* $\psi_n$ *on* $\mathbb{R}$ *by*

$$\psi_n(\lambda) = \frac{S_T(2\omega_n/T, \lambda)}{\|S_T(2\omega_n/T, \cdot)\|_\mu},$$

*and put*

$$
\begin{aligned}
\sigma^{-2}(\omega_n) &= S_T\left(\frac{2\omega_n}{T}, \frac{2\omega_n}{T}\right) \\
(7.2) &= \begin{cases} (2-2H)\Gamma^2(1-H)\left(\dfrac{\omega_n}{2}\right)^{2H} J_{-H}^2(\omega_n)V_T, & \omega_n \neq 0, \\ V_T, & \omega_n = 0. \end{cases}
\end{aligned}
$$

*The functions* $\psi_n$ *form an orthonormal basis of* $\mathcal{L}_T$ *and every function* $\psi \in \mathcal{L}_T$ *can be expanded as*

$$\psi(\lambda) = \sum_{n \in \mathbb{Z}} \sigma(\omega_n)\psi\left(\frac{2\omega_n}{T}\right)\psi_n(\lambda),$$

*the convergence taking place in* $L^2(\mu)$.

PROOF. By the remarks in the beginning of the section, the functions $\psi_n$ are orthogonal, and they have unit norm by construction. Let us prove that the system is complete.

By Corollary 4.2, it suffices to show that every kernel function $\lambda \mapsto S_T(\omega, \lambda)$ is in the closure of the linear span of the $\psi_n$. We claim that, for $\omega, \lambda \in \mathbb{R}$,

$$S_T(2\omega, 2\lambda) = \sum_{n \in \mathbb{Z}} \overline{\psi_n(2\omega)}\psi_n(2\lambda).$$

To prove this, note that we have

$$
\begin{aligned}
(7.3) \quad \sum_{n \in \mathbb{Z}} \overline{\psi_n(2\omega)}\psi_n(2\lambda) &= \frac{S_T(2\omega, 0)S_T(0, 2\lambda)}{S_T(0,0)} \\
&\quad + \sum_{n \neq 0} \frac{S_T(2\omega, 2\omega_n/T)S_T(2\omega_n/T, 2\lambda)}{S_T(2\omega_n/T, 2\omega_n/T)}.
\end{aligned}
$$

By Corollary 6.2, the sum on the right-hand side equals

$$c\left(\frac{T^2\lambda\omega}{4}\right)^H J_{1-H}(T\lambda)J_{1-H}(T\omega)e^{iT(\lambda-\omega)}\sum_{n\neq 0}\frac{1}{(\omega_n-T\lambda)(\omega_n-T\omega)},$$

where $c=(2-2H)\Gamma^2(1-H)V_T$. By the preceding lemma, the sum in the last display equals

$$\frac{1}{T(\lambda-\omega)}\left(\frac{J_{2-H}(T\lambda)}{J_{1-H}(T\lambda)}-\frac{J_{2-H}(T\omega)}{J_{1-H}(T\omega)}\right).$$

In view of (6.4), multiplication by $J_{1-H}(T\lambda)J_{1-H}(T\omega)$ yields

$$\frac{J_{2-H}(T\lambda)J_{1-H}(T\omega)-J_{2-H}(T\omega)J_{1-H}(T\lambda)}{T(\lambda-\omega)}$$

$$=\frac{J_{-H}(T\omega)J_{1-H}(T\lambda)-J_{-H}(T\lambda)J_{1-H}(T\omega)}{T(\lambda-\omega)}$$

$$-\frac{2-2H}{T^2\lambda\omega}J_{1-H}(T\lambda)J_{1-H}(T\omega).$$

Hence, using Corollary 6.2 again, we see that the sum on the right-hand side of (7.3) equals

$$S_T(2\omega,2\lambda)-c\left(\frac{T^2\lambda\omega}{4}\right)^H e^{iT(\lambda-\omega)}\frac{2-2H}{T^2\lambda\omega}J_{1-H}(T\lambda)J_{1-H}(T\omega).$$

Therefore, our claim follows from the fact that the second term in the last display is equal to

$$\frac{S_T(2\omega,0)S_T(0,2\lambda)}{S_T(0,0)}=\frac{\hat{m}_T(2\lambda)\overline{\hat{m}_T(2\omega)}}{\hat{m}_T(0)}$$

[recall that $S_T(0,\lambda)=\hat{m}_T(\lambda)$ and $S_T(0,0)=V_T=\hat{m}_T(0)$]. To check this, evoke expression (2.2).

So, indeed, the functions $\psi_n$ form a complete, orthonormal system. It follows that every $\psi\in\mathcal{L}_T$ can be written as $\psi=\sum\langle\psi,\psi_n\rangle_\mu\psi_n$. By the reproducing property,

$$\langle\psi,\psi_n\rangle_\mu=\frac{\int_\mathbb{R}\psi(\lambda)\overline{S_T(2\omega_n/T,\lambda)}\mu(d\lambda)}{\|S_T(2\omega_n/T,\cdot)\|_\mu}=\frac{\psi(2\omega_n/T)}{\|S_T(2\omega_n/T,\cdot)\|_\mu}.$$

Another application of the reproducing property yields

$$\left\|S_T\left(\frac{2\omega_n}{T},\cdot\right)\right\|_\mu^2=S_T\left(\frac{2\omega_n}{T},\frac{2\omega_n}{T}\right)=\sigma^{-2}(\omega_n).$$

The explicit expression in (7.2) of the normalizing factor $\sigma(\omega_n)$ follows from the second assertion of Corollary 6.2. □

We remark that, instead of the zeros of $J_{1-H}$, we can also use the zeros of $J_{-H}$ to obtain a second orthonormal basis of $\mathcal{L}_T$. Since the reasoning is completely analogous to the $J_{1-H}$ case, we omit the details and mention only that all the consequent results of this paper can be easily reformulated in terms of the zeros of $J_{-H}$ [like in the case $H = 1/2$ of the Brownian motion, where there exist expansions in terms of the zeros of the sine and the cosine, cf., e.g., Yaglom (1987), Section 26.1].

Using the isometry $U : \mathcal{K}_T \to \mathcal{L}_T$, we can now easily obtain an orthonormal basis of the function space $\mathcal{K}_T = L^2([0,T], V)$ and the corresponding series expansion. Note that in the case $H = 1/2$ it reduces to the Fourier series expansion, given for $f \in L^2[0,T]$ by

$$f(u) = \frac{1}{T} \sum_{n \in \mathbb{Z}} \hat{f}\left(\frac{2n\pi}{T}\right) e^{-(2in\pi/T)u}.$$

COROLLARY 7.3. *Let* $\cdots < \omega_{-1} < \omega_0 = 0 < \omega_1 < \cdots$ *be the real zeros of* $J_{1-H}$ *and let* $\varphi$ *be given by* (3.1). *Then the functions*

$$u \mapsto \frac{\overline{\varphi((2\omega_n/T)u)}}{\|\varphi((2\omega_n/T)\cdot)\|_V}, \qquad n \in \mathbb{Z},$$

*form an orthonormal basis of* $\mathcal{K}_T = L^2([0,T], V)$. *Every* $f \in \mathcal{K}_T$ *can be expanded as*

$$f(u) = \sum_{n \in \mathbb{Z}} \frac{Uf(2\omega_n/T)\overline{\varphi((2\omega_n/T)u)}}{\|\varphi((2\omega_n/T)\cdot)\|_V^2},$$

*the convergence taking place in* $L^2([0,T], V)$.

PROOF. Just note that $S_T(\omega, \cdot)$ is the image under the isometry $U : \mathcal{K}_T \to \mathcal{L}_T$ of the function $u \mapsto \overline{\varphi(u\omega)}$. □

**8. Paley–Wiener expansion.** In this section we use the orthonormal basis of $\mathcal{L}_T$ to obtain a series expansion of the fBm. Paley and Wiener (1934) use a series expansion to introduce the standard Brownian motion, which they call the "fundamental random function" [see also Lévy (1965), Section 13 of the Complement]. They first consider the series

$$\sum_{n \in \mathbb{Z}} e^{int} Z_n, \qquad t \in [0, 2\pi],$$

where the $Z_n$ are i.i.d., complex-valued, standard Gaussian random variables. This series corresponds to white noise, but the series does not converge in the usual sense. So instead they consider its formal integral

$$\sum_{n \in \mathbb{Z}} \frac{e^{int} - 1}{in} Z_n.$$

The latter series is shown to converge almost surely and is taken as the definition for the (complex-valued) Brownian motion.

We can now present the extension to the case $H \neq 1/2$ of this expansion. It is natural to consider a complex-valued fBm. This is a centered, complex-valued Gaussian process $X$ with covariance structure

$$\mathbb{E}X_s\overline{X}_t = \tfrac{1}{2}(s^{2H} + t^{2H} - |s-t|^{2H}).$$

THEOREM 8.1.   Let $\cdots < \omega_{-1} < \omega_0 = 0 < \omega_1 < \cdots$ be the real zeros of $J_{1-H}$ and let $Z_n$, $n \in \mathbb{Z}$, be independent, complex-valued Gaussian random variables with mean zero and variance $\mathbb{E}|Z_n|^2 = \sigma^2(\omega_n)$, where $\sigma^2(\omega_n)$ is given by (7.2). Then, with probability one, the series

$$\sum_{n \in \mathbb{Z}} \frac{e^{2i\omega_n t/T} - 1}{2i\omega_n/T} Z_n$$

converges uniformly in $t \in [0,T]$ and defines a complex-valued fBm with Hurst index $H$.

PROOF.   By Theorem 7.2, we have

$$\hat{\mathbb{1}}_t(\lambda) = \sum_{n \in \mathbb{Z}} \sigma(\omega_n)\hat{\mathbb{1}}_t\left(\frac{2\omega_n}{T}\right)\psi_n(\lambda),$$

where $\psi_n$ is a complete orthonormal system in $\mathcal{L}_T$. It follows that, for $s, t \in [0,T]$,

$$(8.1) \qquad \langle \hat{\mathbb{1}}_s, \hat{\mathbb{1}}_t \rangle_\mu = \sum_{n \in \mathbb{Z}} \sigma^2(\omega_n)\hat{\mathbb{1}}_s\left(\frac{2\omega_n}{T}\right)\overline{\hat{\mathbb{1}}_t\left(\frac{2\omega_n}{T}\right)}.$$

By the spectral representation of the fBm, the left-hand side of the display equals the covariance function of the fBm. Hence, the equality shows that the series in the statement of the theorem converges in mean square sense for every $t \in [0,T]$ and defines a Gaussian process with the same covariance structure as the fBm.

The fact that the series converges uniformly with probability one can be deduced from the Itô–Nisio theorem [cf., e.g., Ledoux and Talagrand (1991)]. See, for instance, the proof of Theorem 4.5 of Dzhaparidze and Van Zanten (2004) for details.   □

For $H = 1/2$, we have $\omega_n = n\pi$ [see (3.10)] and $\mathbb{E}|Z_n|^2 = 1/T$, so, indeed, we recover the classical Paley–Wiener expansion of the ordinary Brownian motion in this case.

The Paley–Wiener theorem for the fBm shows that the (complex-valued) fBm on $[0,1]$ can be viewed as the formal integral of the process

$$\sum e^{2i\omega_n t} Z_n.$$

In view of (7.2) and the fact that $z \mapsto \sqrt{z} J_\nu(z)$ is bounded, the latter can be seen as a random signal in which the weight of the component with frequency $\omega_n$ is (up to a constant) approximately equal to $\omega_n^{1/2-H}$.

The real-valued version of the expansion is as follows.

COROLLARY 8.2. *Let $X$, $(Y_n)_{n\in\mathbb{N}}$ and $(Z_n)_{n\in\mathbb{N}}$ be independent, real-valued Gaussian random variables with mean zero and variance*

$$\mathbb{E}X^2 = \frac{1}{V_T}, \qquad \mathbb{E}Y_n^2 = \mathbb{E}Z_n^2 = \frac{\sigma^2(\omega_n)}{2},$$

*where $\sigma^2(\omega_n)$ is given by (7.2). Then, with probability one, the series*

$$tX + \sum_{n=1}^{\infty} \frac{\sin 2\omega_n t/T}{\omega_n/T} Y_n + \sum_{n=1}^{\infty} \frac{(\cos 2\omega_n t/T - 1)}{\omega_n/T} Z_n$$

*converges uniformly in $t \in [0,T]$ and defines a real-valued fBm with Hurst index $H$.*

PROOF. Note that the terms indexed by $n$ and $-n$ in the sum in equation (8.1) are complex conjugates, hence,

$$\langle \hat{\mathbb{1}}_s, \hat{\mathbb{1}}_t \rangle_\mu = \sigma^2(0)st + 2\sum_{n=1}^{\infty} \sigma^2(\omega_n) \Re\left( \hat{\mathbb{1}}_s\left(\frac{2\omega_n}{T}\right) \overline{\hat{\mathbb{1}}_t\left(\frac{2\omega_n}{T}\right)} \right).$$

Since

$$\Re\left( \hat{\mathbb{1}}_s\left(\frac{2\omega_n}{T}\right) \overline{\hat{\mathbb{1}}_t\left(\frac{2\omega_n}{T}\right)} \right)$$

$$= \frac{\sin 2\omega_n s/T}{2\omega_n/T} \frac{\sin 2i\omega_n t/T}{2\omega_n/T} + \frac{(\cos 2\omega_n s/T - 1)}{2\omega_n/T} \frac{(\cos 2\omega_n t/T - 1)}{2\omega_n/T},$$

this shows that the series in the statement of the corollary converges in mean square sense and that the resulting process has the same covariance structure as the fBm. Uniform convergence with probability one can be argued as in the proof of the theorem. □

In the paper Dzhaparidze and Van Zanten (2005) we proved that a related series expansion of the fBm is rate-optimal in the sense of Kühn and Linde (2002). This means that the rate of uniform convergence is $N^{-H}\sqrt{\log N}$. Exactly the same reasoning as in Dzhaparidze and Van Zanten (2005) can be used to prove the optimality of the Paley–Wiener expansion of Theorem 8.1. The main idea is simply to use the asymptotic properties $\omega_n \sim n\pi$ and $J^2_{-H}(\omega_n) \sim 2/n\pi^2$ for $n \to \infty$ [cf. Watson (1944)] to estimate the size of the terms in the expansion. We omit the details and just give the optimality result.

THEOREM 8.3. *The expansion of the fBm of Theorem 8.1 is rate-optimal. It holds that*

$$\mathbb{E} \sup_{t \in [0,T]} \left| \sum_{|n| > N} \frac{e^{2i\omega_n t/T} - 1}{2i\omega_n/T} Z_n \right| \lesssim N^{-H} \sqrt{\log N}.$$

Let us mention that related optimal series expansions for the so-called odd and even parts of the fBm, and, consequently, also of the fBm itself, can be found in Dzhaparidze and Van Zanten (2004). Compared with the Paley–Wiener expansion of Theorem 8.1, the representation of Dzhaparidze and Van Zanten (2004) has the drawback that it requires both the positive zeros of $J_{1-H}$ and of $J_{-H}$.

Another expansion optimal in the sense of Kühn and Linde (2002) is obtained by Ayache and Taqqu (2005). Their construction is of a completely different type. It involves a doubly indexed array of i.i.d. Gaussian random variables $Z_{j,k}$, weighted by functions $t \mapsto 2^{-jH}(\Psi_H(2^j t - k) - \Psi_H(-k))$, defined in terms of the Fourier transform of an appropriate mother wavelet $\Psi$.

Finally, we want to mention the possibility to extend the expansion results to the fractional Brownian sheet (fBs). This can be achieved by taking suitable tensor products like in the paper of Dzhaparidze and Van Zanten (2005), where our earlier double series expansion [cf. Dzhaparidze and Van Zanten (2004)] is extended to the fBs. On extending the Paley–Wiener expansion of Theorem 8.1, the construction is analogous and even simpler because we now have only one sequence of Bessel zeros. The resulting expansion of the fBs is again rate-optimal, as can be shown by using the asymptotic properties of the Bessel function and its positive zeros.

**Acknowledgment.** Thanks to Michael Lifshits for pointing out some silly miscalculations.

CENTER FOR MATHEMATICS
AND COMPUTER SCIENCE
KRUISLAAN 413
1098 SJ AMSTERDAM
THE NETHERLANDS
E-MAIL: kacha@cwi.nl

DEPARTMENT OF MATHEMATICS
FACULTY OF SCIENCES
VRIJE UNIVERSITEIT AMSTERDAM
DE BOELELAAN 1081A
1081 HV AMSTERDAM
THE NETHERLANDS
E-MAIL: harry@cs.vu.nl
URL: www.math.vu.nl/~harry